\documentclass[12pt]{article}
\usepackage{amsmath}
\usepackage{amssymb}
\usepackage{graphicx}
\usepackage[cp1251]{inputenc}
\usepackage[russian]{babel}
\newcommand{\il}[2]{\int\limits_{#1}^{#2}}

\newcommand{\ph}{\phantom{a}}
\newcommand{\phh}{\phantom{aaa}}

\newcommand{\sist}[2]{\left\{
\begin{array}{l}
{#1}\\
\ph\\
{#2}
\end{array}
\right.}

\begin{document}

\vskip 20pt

MSC 34C10

\vskip 20pt

\centerline{\bf A generalization of Sturm's comparison theorem}

\vskip 20 pt

\centerline{\bf G. A. Grigorian}
\centerline{\it Institute  of Mathematics NAS of Armenia}
\centerline{\it E -mail: mathphys2@instmath.sci.am}
\vskip 20 pt

\noindent
Abstract. The Riccati equation method is used to establish a new comparison theorem for systems of two linear first order ordinary differential equation. This result is based on a, so called, concept of ''null-classes'', and is a generalization of Sturm's comparison theorem.
\vskip 20 pt


\vskip 20 pt

{\bf 1. Introduction.} Let $p_j(t), \ph q_j(t), \ph j=1,2$ be real-valued continuous functions on the interval $[a,b]$. Consider the second order linear ordinary differential equations
$$
(p_j(t) \phi')' + q_j(t) \phi = 0, \phh j =1,2. \eqno (1. 1_j)
$$

{\bf Definition 1.1.} {\it Eq. $(1.1_2)$ is called a Sturm majorant for Eq.  $(1.1_1)$ on $[a,b]$ if
$$
p_1(t) \ge p_2(t) > 0 \hskip 12pt  and \hskip 12pt q_1(t) \le q_2(t), \phh t\in [a,b]. \eqno (1.2)
$$
If in addition $q_1(t) < q_2(t)$ or
$$
p_1(t) > p_2(t) > 0 \hskip 12pt  and \hskip 12pt  q_2(t) \ne 0
$$
for some $t\in [a,b]$, then Eq. $(1.1_2)$ is called a strict Sturm majorant for Eq.  $(1.1_1)$ on $[a,b]$.}
The Sturm's famous comparison theorem states (see [1], p. 334)

{\bf Theorem 1.1 (Sturm).} {\it Let  $(1.1_2)$ be a Sturm majorant for $(1.1_1)$ and let $\phi =\phi_1(t) \ne~ 0$ be a solution of Eq.  $(1.1_1)$ having exactly $n (n\ge 1)$ zeroes $t = t_1 < t_2 < ... < t_n$ on $(t_0, t^0]$. Let $\phi = \phi_2(t) \ne 0$ be a solution of Eq.  $(1.1_2)$ satisfying
$$
\frac{p_1(t) \phi_1'(t)}{\phi_1(t)} \ge \frac{p_2(t) \phi_2'(t)}{\phi_2(t)}  \eqno (1.3)
$$
at $t=t_0$ (the expression on the right [or left] of (1.3) is considered to be $+\infty$ if $\phi_2(t_0) = 0$ [or $\phi_1(t_0) =0$]; in particular, (1.3) holds at $t=t_0$ if $\phi_1(t_0) = 0$). Then $\phi_2(t)$ has at least $n$ zeroes on $(t_0, t_n]$. Furthermore $\phi_2(t)$ has at least  $n$ zeroes on $(t_0, t^0)$ if either the inequality (1.3) holds at $t =t_0$ or $(1.1_2)$ is  a strict Sturm majorant for Eq.  $(1.1_1)$ on $[t_0,t_n].$}

Note that if $p_j(t) \ne 0, \ph j=1,2, \ph t\in[a,b]$, then the equations $(1.1_j), \ph j=1,2$ are equivalent (reducible) to the systems
$$
\sist{\phi' = \frac{1}{p_j(t)} \psi}{\psi' = - q_j(t) \phi, \phh t\in [a,b].} \eqno (1.4_j)
$$
$j=1,2$. They are particular cases of more general systems
$$
\sist{\phi' = f_j(t) \psi,}{\psi' = - g_j(t)\phi, \ph t\in [a,b],} \eqno (1.5_j)
$$
$j=1,2$, where $f_j(t), \ph g_j(t), \ph j=1,2$ are real-valued continuous functions on $[a,b]$.
On the other hand obviously the systems $(1.5_j), \ph j=1,2$ are reducible to the systems $(1.1_j), \ph j=1,2$ respectively if in particular $f_j(t) \ne 0, \ph t\in [a,b], \ph j=1,2.$ In this case we can reformulate Theorem 1.1 for the systems  $(1.5_j), \ph j=1,2$ as follows.

{\bf Definition1.2.} {\it The system $(1.5_2)$ is called a Shturm majorant for the system $(1.5_1)$ on $[a,b]$, if
$$
f_2(t) \ge f_1(t) > 0, \ph \mbox{and} \ph g_2(t) \ge g_1(t), \ph t\in [a,b]. \eqno (1.6)
$$
If in addition $g_2(t) > g_1(t)$ or
$$
f_2(t) > f_1(t) > 0 \ph \mbox{and} \ph g_2(t)\ne 0
$$
for some $t\in [a,b]$, then the system $(1.5_2)$ is called a strict Sturm majorant for the system $(1.5_1)$  on $[a,b]$.}

{\bf Theorem 1.2.} {\it Let the system $(1.5_2)$  be a Sturm majorant for the system $(1.5_1)$ and let $(\phi_1(t), \psi_1(t))$ be a solution of $(1.5_1)$  such that $\phi_1(t)$ has exactly $n \ph (n\ge 1)$ zeroes $t=t_1 < ... < t_n$ on $(t_0,t^0]$. Let $(\phi_2(t), \psi_2(t))$ be a solution of $(1.5_2)$  such that
$$
\frac{\psi_1(t_0)}{\phi_1(t_0)} \ge \frac{\psi_2(t_0)}{\phi_2(t_0)} \eqno (1.7)
$$
(the expression in the right [or left] of (1.7) is considered to be $+\infty$, if $\phi_2(t_0) = 0$ [or  $\phi_2(t_0) = 0$]; in particular (1.7) holds if $\phi_1(t_0) = 0$). Then $\phi_2(t)$ has at least $n$ zeroes on $(t_0, t_n]$. Furthermore $\phi_2(t)$ has at least $n$ zeroes on $(t_0,t^0)$ if either inequality (1.7) holds  or $(1.5_2)$ is a strict Sturm majorant for $(1.5_1)$ on $[t_0,t_n].$}

Here arises the following question.  Does Theorem 1.2 remain  valid if  we weaken the conditions (1.6) up to  the following ones?
$$
f_2(t) \ge f_1(t) \ge 0 \hskip 12pt  and \hskip 12pt g_2(t) \ge g_1(t), \phh t\in [a,b]. \eqno (1.8)
$$
The following example shows that answering on this question is a difficult problem.

{\bf Example 1.1.} {\it Set: $f(t) \equiv \sist{0, \ph t\in [0,\pi),}{\sin^2 t, \ph t\in [\pi, 2 \pi],} \ph g(t) \equiv 1, \ph t\in [0, 2\pi]$. Consider the system
$$
\sist{\phi' = f(t) \psi,}{\psi' = g(t) \phi, \ph t\in [0,2\pi].}
$$
It is not difficult to verify that for the non trivial solution $(\phi(t), \psi(t))$ of this system with $\phi(0) = 0, \ph \psi(0) =1$ the function $\phi(t)$ is identically zero on $[0,\pi]$, whereas for all solutions of this system with $\phi(0)\ne 0$ the function   has at most finitely many zeros.}

Thus, in contrast to the solutions $\phi(t)$  of Eq.  $(1.1_j)$, it may turn out that for some solution $(\phi(t), \psi(t)$) of the system  $(1.5_j)$ the function  $\phi(t)$ (which is associated with a solution $\phi(t)$ of  Eq.  $(1.1_j)$) may have an infinite number (in the given example, the continuum) of zeros on a finite interval of variation of $t$. Despite this, the  posed question can be answered positively if in the formulation of Theorem 1 2 the word \ph "zeros" \ph is replaced by the so-called word of \ph "zero-classes" \ph (the concept null-class is introduced in [2] and has the property that for every solution $(\phi(t), \psi(t))$ of the system  $(1.5_j)$ the function $\phi(t) \ph (\psi(t))$ has no more than finite numbed of null-classes on on $[a,b]$; see below) and the ordering symbol $<$ is replaced by another symbol $\prec$, ordering the null-classes (see below). It turns out that for any solution $(\phi(t), \psi(t)$) of the system  $(1.4_j)$, each zero-class of the function   $\phi(t)$ coincides   with one of its zeros and each zero of $\phi(t)$  is its a zero-class.  Therefore, the solution of the above question with the indicated changes: "zeroes" $\to$ "null-classes', $< \ph \to \ph \prec$\ph  gives a generalization of Sturm's theorem.

In this paper we use the Riccati equation method to obtain a generalization of Sturm's comparison theorem (Theorem 1.1),  which is based on the concept of "null-classes".

{\bf 2. Auxiliary propositions}. Along with the systems $(1.5_j), \ph j=1,2$ consider the Riccati equations
$$
y' + f_j(t) y^2 + g_j(t) = 0, \ph t\in [a,b], \ph j=1,2,  \eqno (2.1_j)
$$
and the differential inequalities
$$
\eta' + f_j(t) \eta^2 + g_j(t) \ge 0, \ph t\in [a,b], \ph j=1,2. \eqno (2.2_j)
$$

{\bf Remark 2.1.} {\it Every solution of Eq. $(2.12)$ on an interval $[t_1,t_2) \ph (\subset [a,b])$ is also a solution of the inequality $(2.2_2)$ on $[t_1,t_2)$,}

{\bf Remark 2.2.} {\it If $f_1(t) \ge 0, \ph t\in [t_1,t_2)$, then for every $\lambda \in (- \infty, + \infty)$ the function $\eta_\lambda(t) \equiv \lambda - \il{t_1}{t} g_1(\tau) d \tau, \ph t\in [t_1,r_2)$ is a solution to the inequality $(2.2_1)$ on $[t_1,t_2)$.}

The following comparison theorem plays a crucial role in the proof of the main result.

{\bf Theorem 2.1.} {\it Let $y_2(t)$ be a solution of Eq. $(2.1_2)$ on $[t_0,\tau_0) \ph (\subset [a,b])$ and let $\eta_1(t)$ and $\eta_2(t)$ be solutions of the inequalities $(2.2_1)$ and $(2.2_2)$ respectively on $[t_0,\tau_0)$, moreover suppose that $y_2(t_0) \le \eta_j(t_0), \ph j=1,2$. In addition let the following conditions be satisfied

$f_1(t) \ge 0, \ph t\in [t_0, \tau_0); \ph  y_{(0)} - y_2(t_0) + \il{t_0}{t} \exp\biggl\{\il{t_0}{\tau} f_1(s) (\eta_1(s) + \eta_2(s)) ds\biggr\}\biggl[(f_2(\tau) - f_1(\tau)) y_2^2(\tau) + g_2(\tau) - g_1(\tau)\biggr] d \tau \ge 0, \ph t\in [t_0,\tau_0),$
for some $y_{(0)} \in [y_2(t_0), \eta_1(t_0)]$.
 Then Eq. $(2.1_1)$ has a solution $y_1(t)$ with the initial condition $y_1(t_0) \ge y_{(0)}$, on $[t_0,\tau_0)$;
moreover $y_1(t) \ge y_2(t), \ph t \in [t_0,\tau_0).$}

See the proof in [3].

Besides of Theorem 2.1 for the proof of the main result  we need also in the following three lemmas.

{\bf Lemma 2.1.} {\it Let $f_1(t) \ge 0, \ph t\in [t_0,\tau_0),$ and let $(t_1,t_2)$ be the maximum existence interval for a solution $y(t)$ of Eq. $(2.1_1)$, where $t_0 < t_1 < t_2 < \tau_0.$, Then
$$
\lim\limits_{t\to t_2 - 0} y(t) = - \infty, \phh \lim\limits_{t\to t_1 +0} = + \infty. \eqno (2.3)
$$
}
 See the proof in [3].

{\bf Remark 2.3}. {\it The first equality of (2.3) remains valid also in the case when $(t_1,t_2)$ is not the maximum existence interval, but $y(t)$ is not continuable to the right from the point $t_2$.}

{\bf Lemma 2.2.} {\it Let $f_1(t) \ge 0, \ph t\in [t_0,\tau)$, and let  $(t_{1k}, t_{2k}) (\subset [t_0,\tau_0))$ be the the maximum existence interval for a solution $y_k(t)$ of Eq. $(2.1_1), \ph k=1,2.$ In addition let $y_1(t_1) > y_2(t_1)$ for some $t_1 \in (t_{11}, t_{21}) \cap (t_{12}, t_{22}).$ Then $t_{11} > t_{12}$ and $t_{21} > t_{22}$.}

See the proof in [3].

{\bf Lemma 2.3.} {\it Let $f_1(t) \ge 0, \ph t\in [t_0,\tau_0)$, let $y_0(t)$ be a solution of Eq. $(2.1_1)$ on $[t_0,\tau_0)$, and let $\eta_0(t)$ be a solution of the inequality $(2.2_1)$ on $[t_0,\tau_0)$; moreover let $y_0(t_0) \le \eta_0(t_0)$. Then $y_0(t) \le \eta_0(t), \ph t\in [t_0, \tau_0)$.}

See the proof in [3].

{\bf 3. Main result}. On the set $2^\mathbb{R}$
of subsets of real numbers $\mathbb{R}$ define the order relation $\prec$, assuming $x \prec y$ if and only if for every $t_x \in x \in 2^\mathbb{R}, \ph t_y\in y \in 2^\mathbb{R}$ the inequality $t_x < t_y$ is valid. Let $(\phi(t), \psi(t))$ be a nontrivial solution of the system $(1.5_1)$. Since $\phi(t)$ is a continuous function on $[a,b]$ thats zeroes form a closed set.

{\bf Definition 3.1.} {\it A connected component of zeroes of the function  $\phi(t)$  of a solution $(\phi(t), \psi(t))$ of the system $(1.5_1)$ is called a null-element of $\phi(t)$.}

Let $z(t)$ be a solution of Eq. $(2.1_1)$ with $z(a) =i$. Then $z(t)$ exists on $[a,b]$ and $y(t) \equiv Im \hskip 2pt z(t) > 0, \ph t\in [a,b]$ and
$$
\phi(t) = \frac{\mu}{\sqrt{y(t)}} \sin \biggl(\il{a}{t} f_1(\tau) y(\tau) d \tau + \theta\biggr), \phh t\in [a,b], \eqno (3.1)
$$
where $\mu$ and $\theta$ are some real constants (see [2]).  Let $n(\phi)$ be a null-element of the  function  $\phi(t)$.  By (3.1)
$$
\il{a}{t} f_1(\tau) y(\tau) d \tau + \theta = \pi k_0, \ph t\in N(\phi), \ph k_0 \in \mathbb{Z}. \eqno (3.2)
$$
Hereafter by $[t_1,t_2]$ we mean the set of all points of $\mathbb{R}$ lying between $t_1$ and $t_2$, including themselves.

{\bf Definition 3.2.} {\it Null-elements $N_1(\phi)$ and $N_2(\phi)$ of the function $\phi(t)$ are called congene-\linebreak rous, if for every $t_j \in N(\phi), \ph j=1,2$ the inequality
$$
\biggl|\il{t_1}{t} f_1(\tau) y(\tau) d \tau \biggr| < \pi, \ph  t\in [t_1,t_2]
$$
is valid.}

It was shown in [2] that the congeniality relation between null-elements $N_1(\phi)$ and $N_2(\phi)$ is an equivalence relation.

{\bf Definition 3.3.} {\it An equivalence class of congenerous null-elements $N(\phi)$ of the function of a solution $(\phi(t), \psi(t))$ of the system $(1.5_1)$ is called a null-class of $\phi(t)$.}

{\bf Remark 3.1.} {\it It follows from (3.1) and (3.2) that if $f_1(t) \ge 0, \ph t\in [a,b]$, then every null-class consists of only one null-element.}

It was shown in [2] that for every solution $(\phi(t), \psi(t))$ of the system $(1.5_1)$ the function $\phi(t)$ has a finite number of null-classes, and these null-classes are linearly ordered by $\prec$.

Let $(\phi_j(t), \psi_j(t))$ be a solution of the system $(1.5_j), \ph j=1,2,$ and let  $N_1(\phi_2) = [t_1,\tau_1] \prec N_2(\phi_2) = [t_2,\tau_2] \prec ... \prec N_n(\phi_2) = [t_n,\tau_n]$ be all of the null-classes of the function $\phi_2(t)$ on $[t_0,t^0], \ph t^0 \in (t_0,\tau_0).$

{\bf Definition 3.4.} {\it The quarter $(\phi_1, \psi_1, f_1, g_1)$ is called a majorant of the quarter \linebreak $(\phi_2, \psi_2, f_2, g_2)$ on $[t_0,t^0]$ if the following conditions are satisfied

1. $\psi_1(t_0)/\phi_1(t_0) \le \psi_2(t_0)/\phi_2(t_0)$ (The expression on the left hand [respectively on the right hand] side of the inequality is considered to be equal $+\infty$  if $\phi_1(t_0)=0$ [respectively, if $\phi_2(t_0)=0$] in particular this is in the case if $\phi_2(t_0) = 0$);

2. $f_1(t) \ge f_2(t) \ge 0, \phh t\in [t_0,t^0];$

3. there exists $\xi_k \in (\tau_k, t_{k+1}) \ph (k=0, ..., n-1)$ such that

($3_1$) $g_1(t) \ge g_2(t), \ph t\in [\tau_k,\xi_k] \ph (k=0, ..., n-1)$

($3_2$) any solutions $\eta_{jk}(t)$ of the inequalities $(2.1_j) \ph (j=1,2)$ on $[\xi_k, t_{k+1}]$ such that $\eta_{jk}(\xi_k) > \psi_2(\xi_k) / \phi_2(\xi_k)$  (such solutions always exist by condition 2 and Remark 2.2) satisfy the inequalities
$$
\il{\xi_k}{t}\exp\biggl\{\il{\xi_k}{\tau} f_2(s)\bigl[\eta_{1,k}(s) + \eta_{2,k}(s)\bigr]d s\biggr\}\biggl[g_1(\tau) - g_2(\tau)\biggr] d \tau \ge 0, \ph t\in [\xi_k, t_{k+1}],\ph  k= 0,..., n-1.
$$
In addition, suppose that $f_2(t) > 0, \ph t\in [t_0, t^0]$ and either the strict inequality takes place in condition 1 or at least one of the following conditions is satisfied:

1'. $f_1(t') > f_2(t')$ and $g_1(t') \ne 0$ for some  $t' \in [t_0, t^0],$

2'. $g_1(t') > g_2(t')$ for some $t' \in \bigcup \limits_{k=1}^{n-1} [\tau_k,\xi_k],$

3'.
$
\il{\xi_k}{t_{k+1}}\exp\biggl\{\il{\xi_k}{\tau} f_2(s)\bigl[\eta_{1k}(s) + \eta_{2k}(s)\bigr]d s\biggr\} \bigl[g_1(\tau) - g_2(\tau)\bigr] d \tau > 0$  for some  $k\in \{0, ..., n-1\}$.

\noindent
Then the quarter $(\phi_1, \psi_1, f_1, g_1)$ is called a strict majorant of the quarter $(\phi_2, \psi_2, f_2, g_2)$ on $[t_0,t^0]$.}

{\bf Theorem 3.1 (main result).}  {\it Let $(\phi_2(t), \psi_2(t))$ be a solution  of the system $(1.5_2)$, let $N_1(\phi_2) = [t_1,\tau_1] \prec ... \prec N_n(\phi_2) = [t_n, \tau_n]$ be all of the null-classes of $\phi_2(t)$ on $[t_0,t^0] \linebreak (t^0 \in (t_0,\tau_0))$ and let $(\phi_1(t), \psi_1(t))$ be a solution of the system $(1.5_1)$. If $(\phi_1, \psi_1, f_1, g_1)$ is a majorant for $(\phi_2, \psi_2, f_2, g_2)$ on $[t_0,t^0]$, then the function $\phi_1(t)$ has at least one null-class in $(\tau_k, \tau_{k+1}]$ for each $k =0, ..., n-1.$ In addition, if $(\phi_1, \psi_1, f_1, g_1)$ is a strict majorant for $(\phi_2, \psi_2, f_2, g_2)$ on $[t_0,t^0]$, then the function $\phi_1(t)$ has at least $n$ null-classes on $(t_0,t^0)$.}

Proof. Suppose that the function $\phi_1(t)$ has no zeroes on $[\tau_{k_0}, t_{k_0+1}]$ for some $k_0$.

 Then the function $\psi_1(t)/\phi_1(t)$ exists on $(\tau_{k_0}, \widetilde{\tau}_{k_0 +1}]$, foe some $\widetilde{\tau}_{k_0 +1} > t_{k_0 + 1}$
  and is a solution of Eq. $(2.1_1)$ there. At first consider the case when $k_0 \ne 0$. In this case $(\tau_{k_0}, t_{k_0 + 1})$ is the maximal existence interval for the solution $y_1(t) \equiv \psi_2(t)/ \phi_2(t)$ of Eq. $(2.1_2)$. Let $\eta_{j, k_0}(t) \ph (j=1,2)$
be solutions of the respective inequalities $(2.2_j)$ on $[\xi_{k_0}, t_{k_0 +1}]$ with the initial conditions $\eta_{j,k_0}(\xi_{k_0}) > y_2(\xi_{k_0}), \ph j=1,2$
(in virtue of condition 1 and Remark 2.2 these solutions exist always). Let $\widetilde{y}_2(t)$ be a solution of Eq. $(2.1_2)$ such that $y_2(\xi_{k_0}) < \widetilde{y}_2(\xi_{k_0}) \le \min\limits_{j=1,2}\{\eta_{j, k_0}(\xi_{k_0})\}$ and let $(\widetilde{t}_{k_0}, \widetilde{t}_{k_0 + 1})$ be the maximum existence interval for $\widetilde{y}_2(t)$. By Lemma 2.2 it follows from the inequality $y_2(\xi_{k_0}) < \widetilde{y}_2(\xi_{k_0})$ that $\tau_{k_0} < \widetilde{t}_{k_0}$ and $t_{k_0 + 1} < \widetilde{t}_{k_0 +1}$. We assume that  $\widetilde{y}_2(\xi_{k_0})$ is close enough to $y_2(\xi_{k_0})$ to ensure that $\widetilde{t}_{k_0} \in (t_{k_0}, \xi_{k_0})$ and $\widetilde{t}_{k_0 +1}\in (t_{k_0}, \widetilde{\tau}_{k_0 + 1})$. Since $\widetilde{t}_{k_0}$ is the left endpoint of the maximum existence interval of $\widetilde{y}_2(t)$ by Lemma 2.1 we have $\widetilde{y}_2(\widetilde{t}_{k_0} + 0) = + \infty$ and since $\widetilde{y}_2(\widetilde{t}_{k_0}) < +\infty$ [because of the inclusion $\widetilde{t}_{k_0} \in (\tau_{k_0},\xi_{k_0}]]$ we have
$$
y_1(\zeta_{k_0}) \le \widetilde{y}_2(\zeta_{k_0}), \eqno (3.5)
$$
for some  $\zeta_{k_0} \in (\widetilde{t}_{k_0}, \xi_{k_0})$.  Let $(\phi_1, \psi_1, f_1, g_1)$ be a majorant for $(\phi_2, \psi_2, f_2, g_2)$. Then by virtue of Theorem 2.1 it follows from (3.5) that
$$
y_1(\xi_{k_0}) \le \widetilde{y}_2(\xi_{k_0}) \eqno (3.6)
$$
and since $\widetilde{y}_2(\xi_{k_0}) \le \min\limits_{j=1,2}\{\eta_{j, k_0}(\xi_{k_0})\}$, we have
$$
y_1(\xi_{k_0}) \le \eta_{1, k_0}(\xi_{k_0}), \phh \widetilde{y}_2(\xi_{k_0}) \le  \eta_{2, k_0}(\xi_{k_0}).
$$
By Lemma 2.3 from here it follows that
$$
y_1(t) \le \eta_{1, k_0}(t), \phh \widetilde{y}_2(t) \le  \eta_{2, k_0}(t) \phh t\in [\xi_{k_0}, \widetilde{t}_{k_0 +1}).
$$
Therefore $\widetilde{y}_2(\widetilde{t}_{k_0 +1} - 0) \ge y_1(\widetilde{t}_{k_0 +1}) > - \infty$. Then by virtue of Lemma 2.1 $(\widetilde{t}_{k_0}, \widetilde{t}_{k_0+1})$ is not the maximum existence interval for $\widetilde{y}_2(t)$. The obtained contradiction shows that $\phi_1(t)$ has at least one zero $l_0$ on $(\tau_{k_0}, t_{k_0 = 1}]$, which
belongs to a null-element of the function $\phi_1(t)$. Hence according to Remark 3.1 $l_0$ belongs to a null-class $N(\phi_1)$ of the function $\phi_1(t)$. By (3.1) from condition 2 it follows that if $N(\phi_1) \cap N_{k_0}(\phi_2) \neq \emptyset.$  Then $N(\phi_1) \subset N_{k_0}(\phi_2)$, and, therefore, $N(\phi_1) \subset (\tau_{k_0}, \tau_{k_0 + 1}]$. Now consider the case when $k_0 = 0$. If $\phi_2(t_0) = 0$, then the proof of a null-class of $\phi_1(t)$ in
$(t_0,t_1]$ can be proved by analogy with the proof in the preceding case. Suppose $\phi_2(t_0) \ne 0$. Then by condition 1 we have also $\phi_1(t_0) \neq 0$. Show that $\phi_1(t)$ has at least one zero on $(t_0,t_1]$. Suppose $\phi_1(t) \ne 0, \ph t\in (t_0,t_1]$.  Then $y_1(t)\equiv \psi_1(t)/ \phi_1(t)$ is defined at least on $[t_0,t_1]$ and is a solution of Eq. $(2.1_1)$ on that interval. The function $y_2(t) \equiv \psi_2(t)/\phi_2(t)$ is a solution of Eq. $(2.1_2)$ on $[t_0,t_1)$, in virtue of Remark~ 2.3 $y_2(t_1 - 0) = -\infty$. But on the other hand since $y_2(t_0) \ge y_1(t_0)$ and $(\phi_1, \psi_1, f_1, g_1)$ is a majorant for  $(\phi_2, \psi_2, f_2, g_2)$ it follows from Theorem 2.1 that $y_2(t_1 - 0) \ge y_1(t_1) > - \infty$. The obtained contradiction shows that $\phi_1(t)$ has at least one zero on $(t_0, t_1]$. Then $\phi_1(t)$ has at least one null-class in $[t_0,\tau_1]$.  The first part of the theorem is proved. The second (last) part of the theorem can be proved by analogy of the proof of the second part of Theorem 4.1 from [3], as far as the strict majorant condition implies the reducibility of the systems $(1.5_j), \ph j=1,2$ to the second order linear ordinary differential equations like $(1.1_j), \ph j=1,2$ respectively. The theorem is proved.

{\bf Definition 3.5.} {\it The system $(1.5_1)$ is called a majorant of the system $(1.5_2)$ on $[t_0, t^0] \linebreak(\subset [a,b])$ if the following conditions are satisfied

1.) $f_1(t) \ge f_2(t) \ge 0, \phh t\in [t_0,t^0],$

2). $g_1(t) \ge g_2(t), \phh t\in [t_0, t^0]$.

In addition suppose that $f_2(t) > 0, \ph t\in [t_0,t^0]$ and at least  one of the following conditions is satisfied

1') $f_1(t') > f_2(t')$ and $g_1(t') \ne 0$ for some $t' \in [t_0, t^0]$;

2') $g_1(t') > g_2(t')$ for some  $t' \in [t_0, t^0]$.

Then the system $(1.5_1)$ is called a strict majorant of the system $(1.5_2)$.}

From Theorem 3.1 we immediately get.

{\bf Corollary 3.1.} {\it Let the system $(1.5_1)$ be a majorant for the system $(1.5_2)$. Let \linebreak $(\phi_2(t), \psi_2(t))$ be a nontrivial solution of the system $(1.5_2)$ and let  $\phi_2(t)$ have exactly $n (\ge 1)$ null-classes $N_1(\phi_2) \prec ... \prec N_n(\phi_2)$ on $[t_0,t^0]$. Let $(\phi_1(t), \psi_1(t))$ be a nontrivial solution of the system $(1.5_1)$ satisfying
$$
\frac{\psi_2(t_0)}{\phi_2(t_0)} \ge \frac{\psi_1(t_0)}{\phi_1(t_0)} \eqno (3.7)
$$
(The expression on the right [or left] of (3.7) is considered to be  $+ \infty$, if $\phi_1(t_0) = 0$ [or $\phi_2(t_0) = 0$]; in particular, (3.7) holds if $\phi_2(t_0)= 0$). Then $\phi_1(t)$ has at least $n$ null-classes in $(t_0, t_n)$, where $N_n(\phi_2) = [t_n, \tau_n]$, if either the strict inequality $(3.7)$ holds or $(1.5_1)$ is a strict majorant for $(1.5_2)$ on $[t_0, t_n]$.}

$\phantom{aaaaaaaaaaaaaaaaaaaaaaaaaaaaaaaaaaaaaaaaaaaaaaaaaaaaaaaaaaaaaaaaaaaa} \Box$

{\bf Remark 3.2.} {\it Corollary 3.1 is a generalization of Theorem 1.1.}

\vskip 20 pt

\centerline{ References}

\vskip 15 pt

\noindent
1. Ph. Hartman, Ordinary differential equations, SIAM - Society for industrial and\linebreak \phantom{aaa} applied Mathematics, Classics in Applied Mathematics 38, Philadelphia 2002.

  \noindent
2.  G. A. Grigorian. Oscillatory criteria for the systems of two first - order Linear\linebreak \phantom{a} ordinary differential equations. Rocky Mount. J. Math., vol. 47,\linebreak \phantom{a} Num. 5, 2017, pp. 1497 - 1524

\noindent
3. G. A. Grigorian,  On two comparison tests for second-order linear  ordinary\linebreak \phantom{aa} differential equations (Russian) Differ. Uravn. 47 (2011), no. 9, 1225 - 1240; trans-\linebreak \phantom{aa} lation in Differ. Equ. 47 (2011), no. 9 1237 - 1252, 34C10.

\end{document}